\documentclass[10pt]{article}
\usepackage[all]{xy}
\usepackage{amssymb}

\newenvironment{prf}[1]{\trivlist 
\item[\hskip \labelsep{\bf 
#1.\hspace*{.3em}}]}{~\hspace{\fill}~$\square$\endtrivlist} 
\newenvironment{proof}{\begin{prf}{Proof}}{\end{prf}} 

\newtheorem{theorem}{Theorem}[section]
\newtheorem{lemma}[theorem]{Lemma}
\newtheorem{proposition}[theorem]{Proposition}
\newtheorem{corollary}[theorem]{Corollary}
\newtheorem{remark}[theorem]{Remark}
\newtheorem{definition}[theorem]{Definition}

\newcommand{\Iemph}[1]{\emph{#1}\index{#1}}

\newcommand{\ab}[3]{#1:#2\rightarrow#3}

\newcommand{\hsp}{\hspace{5pt}}
\newcommand{\spec}[1]{ \mathrm{Spec} (#1)}

\newcommand{\xz}{\mathbb{Z}}

\newcommand{\xg}{\mathbb{G}}

\newcommand{\pol}{\mathcal{L}}

\begin{document}

\title{Canonical coordinates on the canonical lift
\thanks{This work was partially supported by the
Australian Research Council grant DP0453134 and by
the Spinoza grant of H.W. Lenstra.}}

\author{Robert Carls \\ \\
School of Mathematics and Statistics \\
University of Sydney \\
NSW 2006 Australia
\\ \\
tel: +61293515775 \\
fax: +61293514534 \\
email: \tt{carls@maths.usyd.edu.au} }

\date{\today}

\maketitle

\begin{abstract}
\noindent
In this article we prove the existence of a canonical theta
structure for the canonical lift of an ordinary abelian variety.
\end{abstract}

\section{Introduction}

The aim of this article is to provide a theoretical basis for
the study of the theta null points of canonical lifts.
We prove that there exists a \emph{canonical theta
structure} for the canonical lift of an ordinary abelian variety.
The canonical theta null point of a canonical lift, which is defined in terms
of the canonical theta structure,
forms an arithmetic invariant.
\newline\indent
In later work we will show how the results of the present article can
be used to compute equations satisfied by the theta
null point of a canonical lift over a $3$-adic ring. Together with
D. Kohel and D. Lubicz we have implemented an algorithm, based on the
latter equations, for the construction of hyperelliptic curves of genus
$2$ with CM over number fields.
\newline\indent
The main result, Theorem \ref{excantheta}, was proven during the last two
years of the author's postgraduate studies at the universities of Leiden
and Groningen.
It extends the result presented in the author's
PhD thesis \cite[Th. 4.1.1]{ca04} to the case of residue
field characteristic $2$.

\section{The canonical theta structure}
\label{thecanonicalthetastructure}

\noindent
For general remarks about the notation see Section \ref{notation}.
Let $R$ be a complete noetherian local ring with
residue field $k$ of characteristic $p>0$
and let $A$ be an abelian scheme over $R$ of relative dimension $g$
having ordinary reduction.
Let $\pol$ be an ample line bundle on $A$.
Let $j \geq 1$ and $q=p^j$.
It is known that there
exists an isogeny of abelian schemes $F:A \rightarrow A^{(q)}$
which is uniquely determined up to unique isomorphism by the condition that
it lifts the relative $q$-Frobenius on the special fibre
(compare \cite[Prop. 2.2.1]{ca04}).
In Section \ref{first} we prove the existence of
an ample line bundle $\pol^{(q)}$ on $A^{(q)}$
satisfying $F^* \pol^{(q)} \cong \pol^{\otimes q}$
which is uniquely determined up to isomorphism by
the condition that
$\pol^{(q)}$ restricted to the special fibre is the $q$-Frobenius
twist of $\pol$.
For a precise statement see Theorem \ref{key1}.
Assume that we are given an isomorphism
\begin{eqnarray}
\label{trivp}
(\xz / q \xz)_R^g \stackrel{\sim}{\rightarrow} A[q]^{\mathrm{et}}
\end{eqnarray}
where $A[q]^{\mathrm{et}}$ denotes the maximal \'{e}tale
quotient of $A[q]$.
\begin{theorem}
\label{excantheta}
Let $\pol$ be an ample symmetric line bundle of degree $1$ on $A$.
There exists a
canonical theta structure of type $(\xz / q \xz )_R^g$
for the pair
\[
\left( A^{(q\delta)}, \big( \pol^{(q\delta)}
\big)^{\otimes q} \right)
\quad \mbox{where} \quad
\delta= \left\{ \begin{array}{l@{, \quad}l}
2 & p=2 \\
1 & p>2
\end{array} \right. 
\]
depending on the isomorphism \rm (\ref{trivp}) \it.
\end{theorem}
\noindent
Theorem \ref{excantheta} will be proven in Section \ref{excanthetaproof}.
For the definition of a theta structure we refer to Section \ref{tstructure}.
For the following statement we assume
that $k$ is perfect and $R$ admits an automorphism lifting the
$p$-th power Frobenius automorphism of $k$.
\begin{corollary}
\label{canlifttheta}
Now assume that $q>2$.
Let $A$ be a canonical lift and
let $\pol$ be an ample symmetric line bundle of degree $1$ on $A$.
There exists a canonical theta structure of type $(\xz / q \xz)^g_R$
for the pair $\big( A, \pol^{\otimes q} \big)$
depending on the isomorphism \rm{(\ref{trivp})}.
\end{corollary}
\noindent
Corollary \ref{canlifttheta}
will be proven in Section \ref{canliftproof}.
The above corollary is expected to hold
without the assumption that $R$
admits a lift of the $p$-th power Frobenius automorphism of $k$.
For a discussion of canonical lifting see \cite{me72},
\cite[Ch. III]{mo95} or \cite[Ch. 2]{ca04}. 

\section{Notation}
\label{notation}

\noindent
Let $R$ be a ring, $X$ an $R$-scheme and $S$ an $R$-algebra.
By $X_{S}$ we denote the base extended scheme $X
\times_{\spec{R}} \spec{S}$.
Let $\mathcal{M}$ be a sheaf on $X$. Then we denote by $\mathcal{M}_S$
the sheaf 
that one gets by pulling back via the projection $X_S \rightarrow X$.
Let $\ab{I}{X}{Y}$ be a morphism of $R$-schemes. Then $I_S$ denotes the
morphism that is induced by $I$ via base extension with $S$.
We use the same symbol for a scheme and the
fppf-sheaf represented by it. By a \emph{group} we mean a group
object in the category of fppf-sheaves.
If a representing object has the property of
being  finite (flat, \'{e}tale, connected, etc.)
then we simply say that it is a finite (flat, \'{e}tale,
connected, etc.) group.
Similarly we will say that
a morphism of groups is finite (faithfully flat, smooth, etc.)
if the groups are representable and the induced morphism of schemes has the
corresponding property.
\newline\indent
A group (morphism of groups) is called finite locally free if
it is finite flat and of finite presentation.
The Cartier dual of a finite locally free commutative
group $G$ will be denoted by $G^D$.
The multiplication by an integer $n \in \xz$ on $G$ will be denoted by
$[n]$.
A finite locally free and surjective morphism between groups is called
an \emph{isogeny}.
By an \Iemph{elliptic curve} we mean an abelian scheme of relative
dimension $1$.
We use the notion of a \Iemph{torsor} in the sense of
\cite[Ch. III, $\S$4, Def. 1.3]{dg70}. 
We only consider torsors for the fppf-topology.

\section{Theta groups and theta structures}
\label{theory}

\noindent
In the following sections we recall some well-known facts
about theta groups and theta structures. We refer to
\cite{mu66}, \cite[Ch. IV, $\S$23]{mu70}, \cite[Ch. 8]{gmav}
and \cite[Ch. V]{mb85} for further details.
Let $R$ be a ring and $G$ a group over $R$.
\begin{definition}
\label{thetagroup}
Assume that there exists a central exact sequence of groups
\[
0 \rightarrow \xg_{m,R} \rightarrow G \rightarrow
H \rightarrow 0,
\]
where $H$ is a commutative finite locally free group whose rank is the
square of an integer. 
Then the group $G$ is called a \Iemph{theta group} over $H$.
\end{definition}
\noindent
By the term \Iemph{central exact sequence} we mean that $\xg_{m,R}$
is mapped into the centre of $G$.
Now let $G$ be a theta group over $H$.
By definition $G$ is a $\xg_{m,R}$-torsor over $H$.
It follows by descent that the group $G$ is representable
by an affine faithfully flat group scheme of finite presentation over $H$
(see \cite[Ch. III, $\S$4, Prop. 1.9]{dg70}).
Let $S$ be an $R$-algebra.
One defines the \Iemph{commutator pairing}
$\mathrm{e}:H \times_R H \rightarrow \xg_{m,R}$
by lifting $x$ and $y$ in $H(S)$ to $\tilde{x}$ and $\tilde{y}$
in $G(S')$, where $S \rightarrow S'$ is a suitable fppf-extension,
and by setting
\[
\mathrm{e}(x,y)= \tilde{x}\tilde{y}\tilde{x}^{-1}\tilde{y}^{-1}.
\]
Because $H$ is abelian we have $\mathrm{e}(x,y) \in \xg_{m,R}(S')$.
Since $\mathrm{e}(x,y)$ does not depend on the choice of $\tilde{x}$
and $\tilde{y}$ it follows by descent that $\mathrm{e}(x,y) \in \xg_{m,R}(S)$.

\subsection{The theta group of an ample line bundle} 
\label{thetalb}

Let $A$ be an abelian scheme over a ring $R$ and $\pol$ a
line bundle on $A$.
Consider the morphism
\[
\varphi_{\pol}:A \rightarrow \mathrm{Pic}^0_{A/R},
\hsp x \mapsto \langle T_x^* \pol \otimes \pol^{-1} \rangle
\]
where $\langle \cdot \rangle$ denotes the class in $\mathrm{Pic}^0_{A/R}$.
We set $\check{A}=\mathrm{Pic}^0_{A/R}$. Note that $\check{A}$
is the dual of $A$ in the category of abelian schemes.
It is well-known that the relative Picard functor of $A$
is representable by an
algebraic space (compare \cite[Ch. 8, Th. 1]{bl90}).
The representability of $\mbox{Pic}^0_{A/R}$ by a scheme
follows from a theorem of M. Raynaud which states that
the categories of abelian
algebraic spaces and abelian schemes coincide
(see \cite[Ch. I, Th. 1.9]{fc90}).
We denote the kernel of the morphism $\varphi_{\pol}$ by $H(\pol)$.
A line bundle $\pol$ on $A$ satisfies $H(\pol)=A$
if and only if its class is in $\mbox{Pic}^0_{A/R}(R)$.
Also it is well-known that if $\pol$ is relatively ample then
$\varphi_{\pol}$ is an isogeny. In the latter case we say that
$\pol$ has degree $d$ if $\varphi_{\pol}$ is fibre-wise of
degree $d$.
Let $S$ be an $R$-algebra.
We define
\[
G(\pol)(S)= \left\{ \hsp (x, \varphi) \hsp | \hsp x \in H(\pol)(S),
\hsp \varphi:\pol_S \stackrel{\sim}{\rightarrow}
T^*_x \pol_S \hsp \right\}.
\]
The functor $G(\pol)$ has the structure of a group given by the group law
\[
\big( (y,\psi),(x,\varphi) \big) \mapsto (x+y, T_x^* \psi \circ \varphi).
\]
There are natural morphisms
\[
G(\pol) \rightarrow H(\pol), \hsp (x,\varphi) \mapsto x \quad \mbox{and}
\quad \xg_{m,R} \rightarrow G(\pol), \hsp \alpha \mapsto (0_A,\tau_{\alpha})
\]
where $0_A$ denotes the zero section of $A$ and $\tau_{\alpha}$ denotes
the automorphism of $\pol$ given by the multiplication with $\alpha$.
The induced sequence of groups
\begin{eqnarray}
\label{thetaseq}
0 \rightarrow \xg_{m,R} \rightarrow G(\pol) \stackrel{\pi}{\rightarrow} H(\pol)
\rightarrow 0
\end{eqnarray}
is central and exact.
Now let $\pol$ be relatively ample of degree $d$.
Then $H(\pol)$ is finite locally
free of order $d^2$
and hence $G(\pol)$ is a theta group. The commutator pairing on $H(\pol)$
as defined above will be denoted by $\mathrm{e}_{\pol}$.
One can show that the pairing $\mathrm{e}_{\pol}$ is perfect.
The perfectness is
equivalent to the fact that the centre of $G(\pol)$ equals $\xg_{m,R}$.

\subsection{Descent of line bundles along isogenies}
\label{descisog}

Let $R$ be a ring.
Let $\ab{I}{A}{B}$ be an isogeny of abelian schemes over $R$
and $K$ its kernel.
Assume we are given a relatively ample line bundle $\pol$ on $A$
and $K \subseteq H(\pol)$.
Define $G'$ by the commutative diagram
\begin{eqnarray}
\label{belang}
\xymatrix{ 0 \ar@{->}[r] & \xg_{m,R} \ar@{<-}[d]^{\mathrm{id}}
\ar@{->}[r] & G(\pol) \ar@{->}[r] \ar@{<-}[d] &
H(\pol) \ar@{->}[r] \ar@{<-}[d]_i & 0 \\
0 \ar@{->}[r] & \xg_{m,R} \ar@{->}[r] &
G' \ar@{->}[r]^{\pi} & K \ar@{->}[r] & 0,}
\end{eqnarray}
where the second row is the pull back of the first via the inclusion
$K \stackrel{i}{\hookrightarrow} H(\pol)$, i.e. the right hand square is
Cartesian.
Let $U$ be an $R$-algebra and
$\mathcal{M}$ a line bundle on $B_U$.
Suppose we are given an
isomorphism $\alpha:I_U^* \mathcal{M} \stackrel{\sim}{\rightarrow}
\pol_U$.
We define a morphism $s_{\alpha}:K_U \rightarrow G_U'$
by mapping $x \in K(W)$, where $W$ is a $U$-algebra,
to $\big( x , T_x^*\alpha_{W} \circ \alpha_{W}^{-1} \big)$.
This is well-defined because $T^*_xI_{W}^* \mathcal{M}_{W}
=I_{W}^* \mathcal{M}_{W}$.
It is clear that $\pi_U \circ s_{\alpha} = \mathrm{id}$ where $\pi$
is as in diagram (\ref{belang}).
We define
\[
S_K(U)
=\{ \hsp s:K_U \rightarrow G_U' \hsp | \hsp \pi_U \circ s = \mathrm{id} \hsp \}
\]
and denote by $D_\pol(U)$
the set of isomorphism classes of line bundles $\mathcal{M}$
on $B_U$ such that $I_U^* \mathcal{M} \cong \pol_U$.
The following classical result about the descent of line bundles
was proven by Alexander Grothendieck.
\begin{proposition}
\label{groth}
The functorial map
\[
S_K(U) \rightarrow D_{\pol}(U), \hsp (\mathcal{M}, \alpha) \mapsto s_{\alpha}
\]
establishes an isomorphism of functors $S_K \stackrel{\sim}{\rightarrow}
D_{\pol}$.
\end{proposition}
\noindent
Compare \cite[ Ch. IV, $\S$23, Th. 2]{mu70} or \cite[Ch. 6.1,
Th. 4]{bl90}.

\subsection{Theta structures}
\label{tstructure}

In the following we define the standard theta group of a given type.
Let $K$ be a commutative finite locally free group of square
order over a base ring $R$.
We set
$H(K)=K \times_R K^D$
and define a group law on $G(K)=\xg_{m,R} \times_R
H(K)$ by setting
\[
( \alpha_1, x_1, l_1 ) \ast ( \alpha_2, x_2, l_2)=
(\alpha_1 \cdot \alpha_2 \cdot l_2(x_1),x_1 + x_2,l_1 \cdot l_2).
\]
We have an exact sequence of groups
\[
0 \rightarrow \xg_{m,R} \rightarrow G(K) \rightarrow H(K) \rightarrow 0
\]
where the left hand map is given by $\alpha \mapsto (\alpha,0,1)$
and the right hand map is the projection on $H(K)$.
The centre of $G(K)$ is given by $\xg_{m,R}$.
We conclude that $G(K)$ is a theta group.
We denote the corresponding commutator pairing by $\mathrm{e}_K$. 
Using the definition of the multiplication in $G(K)$
one computes
\begin{eqnarray}
\label{explicitpair}
\mathrm{e}_K \big( (x_1,l_1), (x_2,l_2) \big)=\frac{l_2(x_1)}{l_1(x_2)}.
\end{eqnarray}
We remark that $\mathrm{e}_K$ is a perfect pairing.
Now assume we are given an abelian scheme $A$ over $R$
and a relatively ample line bundle $\pol$ on $A$.
\begin{definition}
A \Iemph{theta structure} of type $K$ for the pair
$(A,\pol)$ is an isomorphism $\Theta:G(K) \stackrel{\sim}{\rightarrow} G(\pol)$
making the diagram
\[
\xymatrix{
\xg_{m,S} \ar@{->}[r] \ar@{<-}[d]^{\mathrm{id}}
& G(\pol) \ar@{<-}[d]^{\Theta} \\
\xg_{m,S} \ar@{->}[r] & G(K) }
\]
commutative. Here the horizontal arrows are the natural inclusions.
\end{definition}
\noindent
Next we want to give another characterisation of a theta structure.
\begin{definition}
A \Iemph{Lagrangian decomposition} for $H(\pol)$ of type $K$
is an isomorphism
\[
\delta: H(K) \stackrel{\sim}{\rightarrow} H(\pol),
\]
which is compatible with the commutator pairings $\mathrm{e}_{\pol}$
and $\mathrm{e}_K$.
\end{definition}
\noindent
Let $\delta$ be a Lagrangian decomposition for $H(\pol)$ of
type $K$.
We can consider $K$ and $K^D$ as subgroups of $H(\pol)$
via $\delta$.
Assume we are given a pair $(u,v)$ where $u$ and $v$ are sections
of the pull back of the extension
\begin{eqnarray}
\label{thetaextseq}
0 \rightarrow \xg_{m,R} \rightarrow G(\pol)
\stackrel{\pi}{\rightarrow} H(\pol)
\rightarrow 0
\end{eqnarray}
along the inclusions $K \hookrightarrow H(\pol)$ and $K^D
\hookrightarrow H(\pol)$, respectively.
We define a morphism
$\Theta_{u,v}:G(K) \rightarrow G(\pol)$ by
$\Theta_{u,v}(\alpha,x,l)= \alpha \cdot v(l) \cdot u(x)$.
\begin{proposition}
\label{newtheta}
The map
\begin{eqnarray}
\label{bijectie}
( \delta, u, v ) \mapsto \Theta_{u,v}
\end{eqnarray}
gives a bijection between the set of triples as above and the
set of theta structures for $(A,\pol)$ of type $K$.
\end{proposition}
\begin{proof}
First we have to show that the map (\ref{bijectie}) is well-defined.
We claim that $\Theta_{u,v}$ is a theta structure of type
$K$ for $(A,\pol)$.
We have
\[
\Theta_{u,v} \big( (\alpha_1,x_1,l_1) \ast (\alpha_2, x_2,
  l_2) \big) = \alpha_1 \cdot \alpha_2 \cdot l_2(x_1) \cdot v(l_1) \cdot v(l_2)
\cdot u(x_1) \cdot u(x_2).
\]
By the definition of the pairing $\mathrm{e}_{\pol}$ it follows that  
\[
v(l_2) \cdot u(x_1) =\mathrm{e}_{\pol} \big( \delta(l_2),\delta(x_1) \big) \cdot
u(x_1) \cdot v(l_2).
\]
Since $\delta$ is a Lagrangian decomposition we have
\[
\mathrm{e}_{\pol} \big( \delta(0,l_2),\delta(x_1,1) \big)= \mathrm{e}_K \big( (0,l_2),
(x_1,1) \big)
=\frac{1}{l_2(x_1)}.
\]
The right hand equality follows by (\ref{explicitpair}).
This proves that $\Theta_{u,v}$ is a morphism of groups.
Clearly $\Theta_{u,v}$ is $\xg_{m,R}$-equivariant.
\newline\indent
Next we prove that $\Theta_{u,v}$ is an isomorphism by giving an
inverse.
Let $g$ be a point of $G(\pol)$. Then we have $\pi(g)=\delta(x_g,l_g)$ for
uniquely determined
$x_g \in K$ and $l_g \in K^D$. Here $\pi$ denotes the projection map
of the extension (\ref{thetaextseq}).
Now $g$ and $\Theta_{u,v}(1,x_g,l_g)$ both lift $\delta(x_g,l_g)$.
Hence they differ by a unique scalar $\alpha_g$, i.e.
$g= \Theta_{u,v}(\alpha_g,x_g,l_g)$.
An inverse of $\Theta_{u,v}$ is given by the morphism
$g \mapsto (\alpha_g, x_g, l_g)$.
\newline\indent
In order to complete the proof of
Proposition \ref{newtheta} it is sufficient to
give an inverse of the map (\ref{bijectie}).
Assume we are given a theta structure $\Theta$ of type $K$
for the pair
$(A,\pol)$.
The isomorphism $\Theta$ induces an isomorphism $\delta_{\Theta}:H(K) \stackrel{\sim}
{\rightarrow} H(\pol)$. By the definition of the commutator pairing
it follows that the isomorphism $\delta_{\Theta}$ is a Lagrangian decomposition.
There are two natural sections of the natural projection $G(\pol) \rightarrow H(\pol)$
over $K$ and $K^D$ given by
\[
u_{\Theta}:(x,1) \mapsto \Theta \big( 1,x,1 \big) \quad \mbox{and}
\quad v_{\Theta}:(0,l) \mapsto \Theta \big( 1,0,l \big),
\]
respectively.
Here we consider $K$ and $K^D$ as subgroups of $H(\pol)$
via $\delta_{\Theta}$. An inverse of (\ref{bijectie}) is given by
$\Theta \mapsto \big( \delta_{\Theta},u_{\Theta},v_{\Theta} \big)$.
This finishes the proof of the proposition.
\end{proof}

\section{Descent along lifts
of relative Frobenius and Verschiebung}
\label{statements}

\noindent
In the following we recall some facts about
the existence of Frobenius lifts and the descent of line bundles
along lifts of Frobenius and Verschiebung.
Theorem \ref{key1} and \ref{key2} are known to the experts
but they are not yet available in the literature.
We prove them in Section \ref{first} and \ref{second}.
\newline\indent
Let $R$ be a complete noetherian local ring with
residue class field $k$ of characteristic $p>0$
and $A$ an abelian scheme
having ordinary reduction.
Let $j \geq 1$ and $q=p^j$.
It is known that there exists an abelian scheme $A^{(q)}$ over $R$
and a commutative diagram of isogenies
\begin{eqnarray*}
\xymatrix{
  A \ar@{->}[r]^F \ar@{->}[d]_{[q]} & A^{(q)} \\
  A \ar@{<-}[ur]_V }
\end{eqnarray*}
such that
$F_k$ equals the relative $q$-Frobenius
(compare \cite[Prop. 2.2.1]{ca04}).
The latter condition determines $F$ uniquely.
The kernel of $F$ is given
by $A[q]^{\rm{loc}}$ which is defined to be the connected
component of $A[q]$.
The condition that $F_k$ equals the relative $q$-Frobenius means that
there exists a commutative diagram
\[
\xymatrix{ A_k \ar@/^/[rrd]^{f_q}
  \ar@/_/[rdd] \ar@{->}[rd]|-{F_k} & & \\
  & \big( A^{(q)} \big)_k \ar@{->}[r]^{\mathrm{pr}} \ar@{->}[d] & A_k
  \ar@{->}[d] \\
  & \spec{k} \ar@{->}[r]^{f_q} & \spec{k} }
\]
where $f_q$ denotes the absolute $q$-Frobenius, the vertical maps are
the structure maps and
the square is Cartesian.
Let $\mathcal{L}$ be a line bundle on $A$.
We have a natural isomorphism
\begin{eqnarray}
\label{natsect}
F^*_k  \mathrm{pr}^* \pol_k = f_q^* \pol_k \stackrel{\sim}{\rightarrow}
\pol_k^{\otimes q}
\end{eqnarray}
given by $l \otimes 1 \mapsto l^{\otimes q}$.
\begin{theorem}
\label{key1}
Assume that $\pol$ is an ample line bundle on $A$.
There exists a line bundle $\pol^{(q)}$ on
$A^{(q)}$ determined uniquely up to isomorphism by
the following two conditions:
\[
\left( \pol^{(q)} \right)_k \cong \mathrm{pr}^* \pol_k \quad
\mbox{and} \quad F^* \pol^{(q)} \cong \pol^{\otimes q}.
\]
Moreover, the line bundle $\pol^{(q)}$ is ample and has the same degree
as $\pol$.
\end{theorem}
\noindent
A proof of Theorem \ref{key1} is presented in Section \ref{first}.
\begin{theorem}
\label{key2}
Assume that $\pol$ is an ample symmetric line bundle on $A$.
\begin{enumerate}
\item
Let $p>2$. There exists an isomorphism
\begin{eqnarray}
\label{keyiso}
V^* \pol \stackrel{\sim}{\rightarrow}
\big( \pol^{(q)} \big)^{\otimes q}.
\end{eqnarray}
\item
\label{casep2}
Let $p=2$. Assume we are given an isomorphism
\begin{eqnarray}
\label{splittingg}
A[2] \stackrel{\sim}{\rightarrow} A[2]^{\mathrm{loc}} \times
A[2]^{\mathrm{et}}.
\end{eqnarray}
There exists a line bundle $\pol_0$ on $A$ with
$\langle \pol_0 \rangle \in \mathrm{Pic}^0_{A/R}[2](R)$
such that
\[
V^* \big( \pol \otimes \pol_0 \big) \stackrel{\sim}{\rightarrow}
\big( \pol^{(q)} \big)^{\otimes q}.
\]
The class of $\pol_0$ depends on the isomorphism
\rm (\ref{splittingg}) \it.
\end{enumerate}
\end{theorem}
\noindent
A proof of Theorem \ref{key2} will be given in Section \ref{second}.
The isomorphism (\ref{keyiso}) does not always exist in the case $p=2$.
This is illustrated by the following example. 
\newline\newline
\bfseries Example: \mdseries
We assume $k$ to be an algebraically closed field of characteristic $2$.
Let $E$ be an ordinary elliptic curve over $k$ and $Q_2$ the unique
non-zero point in $E^{(2)}[2](k)$. Note that $Q_2$ is
a generator of the kernel
of the Verschiebung $V:E^{(2)} \rightarrow E$.
We have
\[
V^*(0_E)=(0_{E^{(2)}})+(Q_2) \not\sim 2 \cdot (0_{E^{(2)}})
\]
where $0_E$ and $0_{E^{(2)}}$ denote the zero sections of $E$
and $E^{(2)}$ and $\sim$ stands for
linear equivalence of Weil divisors.
Let $Q$ be the unique non-zero point in $E[2](k)$
and $R \in E^{(2)}[4](k)$ such that
$2R=Q_2$. Note that $V(R)=Q$. We have
\[
V^* (Q) = (R) + (R+ Q_2) \sim 2 (0_{E^{(2)}}).
\]

\section{The proofs}

\noindent
In the following
we prove the results of Section \ref{thecanonicalthetastructure}
and Section \ref{statements}.

\subsection{Proof of Theorem \ref{key1}}
\label{first}

In the following we prove Theorem \ref{key1}.
We use the notation of Section \ref{statements}.
Let $A$ be an abelian scheme over $R$
having ordinary reduction and let $\pol$ be an ample line bundle on $A$.
Let $K=A[q]^{\mathrm{loc}}$ and
let $G'$ be defined by the Cartesian diagram
\[
\xymatrix{ G(\pol^{\otimes q}) \ar@{->}[r] \ar@{<-}[d] &
  H(\pol^{\otimes q}) \ar@{<-}[d] \\
  G' \ar@{->}[r] & K. }
\]
\begin{remark}
\label{commutativity}
The group $G'$ is commutative.
\end{remark}
\begin{proof}
The commutativity of $G'$ is equivalent to the condition
that the commutator pairing $\mathrm{e}_{\pol^{\otimes q}}:H(\pol^{\otimes q}) \times
H(\pol^{\otimes q}) \rightarrow \xg_m$
is trivial on $K$.
Let $T$ be an $R$-algebra and $x \in K(T)$.
Then the map $e_{\pol^{\otimes q}}(x, \cdot): K_T \rightarrow \xg_{m,T}$
is a $T$-valued point of $K^D \cong \check{A}[q]^{\mathrm{et}}$
where $\check{A}$ denotes the dual abelian scheme.
The morphism $K \rightarrow \check{A}[q]^{\mathrm{et}}$
given by $x \mapsto e_{\pol^{\otimes q}}(x,\cdot)$ is equal to
the zero morphism since the image of $K$
is connected and hence equals
the image of the unit section in $\check{A}[q]^{\mathrm{et}}$ which forms
a connected component.
\end{proof}
\noindent
The main ingredient in the proof of Theorem \ref{key1} is the
following result.
\begin{lemma}
\label{torsoragain}
The functor $S_K$ (defined in Section \ref{descisog}) is a
$K^D$-torsor over $R$.
\end{lemma}
\begin{proof}
Let $U$ be an $R$-algebra.
The group $K^D(U)=\underline{\mbox{Hom}}(
K,\xg_m)(U)$ acts on
$S_K(U)$ by translation.
This action is transitive and faithful in the case that $S_K(U)$ is
non-empty.
Consider the extension
\begin{eqnarray}
\label{overk}
0 \rightarrow \xg_{m,R} \rightarrow G' \rightarrow K
\rightarrow 0.
\end{eqnarray}
It remains to show that (\ref{overk}) has a section over an fppf-extension
of $R$.
Taking $q$-torsion we get an exact sequence
\begin{eqnarray}
\label{ptorsext}
0 \rightarrow \mu_{q,R} \rightarrow G'[q] \rightarrow K \rightarrow 0.
\end{eqnarray}
The exactness of (\ref{ptorsext})
follows from the Snake Lemma and the exactness of the Kummer
sequence.
By Remark \ref{commutativity} the group $G'[q]$ is commutative.
Note that $G'[q]$ is a $\mu_{q,R}$-torsor over $K$. It follows by
\cite[Ch. III, $\S$4, Prop. 1.9]{dg70} that the group $G'[q]$
is finite locally free.
Applying Cartier duality to (\ref{ptorsext}) we get an exact sequence
\begin{eqnarray}
\label{etext}
0 \rightarrow K^D \rightarrow G'[q]^D \stackrel{\pi}{\rightarrow} \xz/ q \xz
\rightarrow 0.
\end{eqnarray}
We remark that
(\ref{ptorsext}) is split if and only if (\ref{etext}) is split.
We can lift a generator of $\xz / q \xz$ to an
element $x \in G'[q]^D(R')$ where $R \rightarrow R'$ is a suitable
fppf-extension.
Clearly $x$ has order $q$.
This gives a splitting of (\ref{etext}) over $R'$. 
As a consequence the extensions (\ref{overk})--(\ref{ptorsext})
are split over $R'$.
\end{proof}
\begin{corollary}
\label{ett}
The functor $S_K$
is representable by a finite \'{e}tale scheme.
\end{corollary}
\begin{proof}
The representability by a finite locally free
$R$-scheme follows from Lemma \ref{torsoragain}
and \cite[Ch. III, $\S$4, Prop. 1.9]{dg70}.
Since $A$ has ordinary reduction, we have $K^D=\underline{\mbox{Hom}}_R(
K,\xg_m) \cong \check{A}[q]^{\mathrm{et}}$ where $\check{A}$ denotes
the dual of $A$.
It follows by descent that $S_K$ is \'{e}tale.
\end{proof}
\noindent
Now we can complete the proof of Theorem \ref{key1}.
We have already seen that there exists
a canonical $k$-rational
point of $S_K$ given by
the isomorphism (\ref{natsect}).
By Corollary \ref{ett} and
the theory of finite \'{e}tale schemes over Henselian local
rings there exists a unique $R$-rational point of
$S_K$ reducing to the above $k$-rational point.
The first part of the
claim of Theorem \ref{key1} now follows from Proposition \ref{groth}.
The second part of the claim states
that the line bundle $\pol^{(q)}$ is ample and of the same
degree as $\pol$.
It suffices to verify the claim
on the special fibre.
It is obvious from the construction that
$\pol^{(q)}_k$ is ample and has the same degree as $\pol_k$.
This finishes the proof of Theorem \ref{key1}.

\subsection{Proof of Theorem \ref{key2}}
\label{second}

We use the notation of Section \ref{statements}.
We set
\[
\mathcal{\pol}'=\big(V^* \pol \big)^{-1} \otimes \big( \pol^{(q)}
\big)^{\otimes q}.
\]
We have $\langle \mathcal{\pol}' \rangle \in \mathrm{Pic}^0_{A^{(q)}/R}(R)$.
In order to prove the proposition we have to show that
$\pol'$ is trivial.
By the symmetry of $\pol$ we conclude that
\begin{eqnarray*}
F^*(V^*\pol) \cong [q]^* \pol \cong \pol^{ \otimes q^2}.
\end{eqnarray*}
Together with Theorem \ref{key1} this implies that
$F^* \pol'$ is trivial on $A$.
This means that $\langle \pol' \rangle$ is in the kernel of the dual
$\check{F}=\mathrm{Pic}^0(F)$ of $F$. The group
$\mathrm{Ker}(\check{F})$
is the Cartier dual of $\mathrm{Ker}(F)$
and hence is annihilated by the isogeny $[q]$.
As a consequence $\langle \pol' \rangle$ has order dividing $q$.
Since we have assumed $\pol$ to be symmetric it follows
that $\pol^{(q)}$ and $\pol'$ are symmetric. We conclude that
\[
\langle \pol' \rangle \in \mathrm{Pic}^0_{A^{(q)}/R}[2](R).
\]
By the above discussion the element $\langle \pol' \rangle$
has order dividing the greatest common divisor of $q$ and $2$.
If $p>2$ we conclude that $\pol'$ is trivial and Theorem \ref{key2}
is proven.
\newline\indent
Assume now that $p=2$. We claim that
there exists a line bundle $\pol_0$ on $A$ with
$\langle \pol_0 \rangle \in \mathrm{Pic}^0_{A/R}[2](R)$ such that
$V^* \pol_0 \cong \pol'$.
We set $\check{A}=\mathrm{Pic}^0_{A/R}$ and $\check{V}=\mathrm{Pic}^0(V)$.
Then $\big( \check{A} \big)^{(q)}
= \mathrm{Pic}^0_{A^{(q)}/R}$.
The isogeny $\check{V}:\check{A} \rightarrow \check{A}^{(q)}$ induces a
morphism of connected-\'{e}tale sequences
\begin{eqnarray}
\label{connetsplit}
\xymatrix{ 0 \ar@{->}[r] & \check{A}[2]^{\mathrm{loc}}
  \ar@{->}[r] \ar@{->}[d]^{\check{V}[2]^{\mathrm{loc}}} &
  \check{A}[2] \ar@{->}[r]^{\pi} \ar@{->}[d]^{\check{V}[2]} &
  \check{A}[2]^{\mathrm{et}}
  \ar@{->}[r] \ar@{->}[d]^{\check{V}[2]^{\mathrm{et}}} & 0 \\
  0 \ar@{->}[r] & \check{A}^{(q)}[2]^{\mathrm{loc}} \ar@{->}[r] &
  \check{A}^{(q)}[2] \ar@{->}[r]^{\pi^{(q)}} &
  \check{A}^{(q)}[2]^{\mathrm{et}}
  \ar@{->}[r] & 0. }
\end{eqnarray}
Note that one cannot embed a non-zero finite \'{e}tale $R$-group into
a connected one.
It follows by the Snake Lemma that $\check{V}[2]^{\mathrm{et}}$
is a monomorphism. Comparing ranks we conclude that
$\check{V}[2]^{\mathrm{et}}$ is an isomorphism. 
By the cokernel property of the morphism $\check{V}[2]^{\mathrm{et}} \circ \pi$
there exists a section
$s:\check{A}^{(q)}[2]^{\mathrm{et}} \rightarrow \check{A}^{(q)}[2]$
of the projection $\pi^{(q)}$.
The image of $\check{A}^{(q)}[2]^{\mathrm{et}}$
under $\pi^{(q)}$ coincides with the kernel of $\check{F}$. Recall
that the latter contains $\langle \pol' \rangle$.
As a consequence there exists an $x \in \check{A}^{(q)}[2]^{\mathrm{et}}$
such that $s(x)$ equals $\langle \pol' \rangle$.
Using the isomorphism (\ref{splittingg}) we map the point 
$(\check{V}[2]^{\mathrm{et}})^{-1}(x)$ to an element of $\check{A}[2]$
whose image under $\check{V}[2]$ equals $\langle \pol' \rangle$.
This proves our claim and completes the proof of the theorem.

\subsection{Proof of Theorem \ref{excantheta}}
\label{excanthetaproof}

We use the notation of Section \ref{thecanonicalthetastructure}.
Let $A$ be an abelian scheme of relative dimension $g$
over $R$ having ordinary reduction and
$\pol$ an ample line bundle of degree $1$ on $A$.
Let $K=(\xz / q \xz)_R^g$.
Assume we are given an isomorphism
\begin{eqnarray}
\label{isooverr}
K \stackrel{\sim}{\rightarrow} A[q]^{\mathrm{et}}.
\end{eqnarray}
The isogeny $F:A \rightarrow A^{(q)}$ induces a
commutative diagram
\begin{eqnarray*}
\label{imp}
\xymatrix{ 0 \ar@{->}[r] & A[q]^{\mathrm{loc}}
  \ar@{->}[r] \ar@{->}[d]^{F[q]^{\mathrm{loc}}} &
  A[q] \ar@{->}[r] \ar@{->}[d]^{F[q]} &
  A[q]^{\mathrm{et}}
  \ar@{->}[r] \ar@{->}[d]^{F[q]^{\mathrm{et}}} & 0 \\
  0 \ar@{->}[r] & A^{(q)}[q]^{\mathrm{loc}} \ar@{->}[r] &
  A^{(q)}[q] \ar@{->}[r] &
  A^{(q)}[q]^{\mathrm{et}}
  \ar@{->}[r] & 0. }
\end{eqnarray*}
The induced morphism $F[q]^{\mathrm{et}}$ is an isomorphism.
Composing the isomorphism (\ref{isooverr}) with $F[q]^{\mathrm{et}}$
we get an isomorphism
$m:K \stackrel{\sim}{\rightarrow} A^{(q)}[q]^{\mathrm{et}}$.
The isomorphism $F[q]^{\mathrm{et}}$ induces a unique section
$r:A^{(q)}[q]^{\mathrm{et}} \rightarrow A^{(q)}[q]$
of the natural projection $A^{(q)}[q] \rightarrow
A^{(q)}[q]^{\mathrm{et}}$.
We define $t=r \circ m$ and
set $H=A^{(q)}[q]$, $C=A^{(q)}[q]^{\mathrm{loc}}$ and
$E=A^{(q)}[q]^{\mathrm{et}}$.
Let $\mathrm{e}(\cdot,\cdot)$ denote the commutator pairing
on
\[
H=H\left( \big( \pol^{(q)} \big)^{\otimes q} \right).
\]
Since $\mathrm{e}$ is a perfect pairing it induces an
isomorphism
$\varphi:H \stackrel{\sim}{\rightarrow} H^D,
\hsp x \mapsto \mathrm{e}(x, \cdot)$.
Note that $C$ is mapped to the connected component of $H^D$.
As a matter of fact the connected component of $H^D$ is given by $E^D$.
Hence the isomorphism $\varphi$ induces isomorphisms
$\alpha:C \stackrel{\sim}
{\rightarrow} E^D$ and
$\beta:E \stackrel{\sim}{\rightarrow} C^D$
on the local and \'{e}tale part of $H$.
We define
$k=-(\alpha^{-1} \circ m^{-D}):
K^D \stackrel{\sim}{\rightarrow} C$
and set $s=i \circ k$.
\begin{lemma}
\label{lagr}
The morphism $\delta=s \oplus t$ is a Lagrangian decomposition of type $K$
for $H$.
\end{lemma}
\begin{proof}
Consider the commutative diagram
\[
\xymatrix{
K^D \ar@{->}[r]^{k} \ar@{->}[d]
\ar@{->}[rd]^{s}
& C \ar@{->}[r]^{\alpha} \ar@{->}[d]^i & E^D
\ar@{->}[r]^{m^D} \ar@{->}[d]^{p^D}
& K^D \ar@{->}[d] \\
K \times K^D \ar@{->}[r]^{\delta} \ar@{->}[d]
& H \ar@{->}[r]^{\varphi} \ar@{->}[d]^p
& H^D \ar@{->}[r]^{\delta^D} \ar@{->}[d]^{i^D}
& K \times K^D \ar@{->}[d] \\
K \ar@{->}[ru]^{t} \ar@{->}[r]^{m}
& E \ar@{->}[r]^{\beta} & C^D \ar@{->}[r]^{k^D}
& K.}
\]
By definition we have
$m^D \circ \alpha \circ k =
-\mathrm{id}$.
Since the pairing $\mathrm{e}$ is alternating it follows that
$\varphi^{D}=-\varphi$.
As a consequence we have
$\beta^D=-\alpha$.
Hence
\[
k^D \circ \beta \circ m = \big( m^D \circ (- \alpha) \circ k \big)^D
=\mathrm{id}.
\]
The commutator pairing $\mathrm{e}_K$ on $K \times K^D$
gives an isomorphism
\[
\tau : K \times K^D \rightarrow K \times K^D, \hsp z \mapsto \mathrm{e}_K \big( z,
\cdot ).
\]
One computes $\tau \big( (x,l) \big) = ( x, l^{-1})$.
We conclude that
$\tau = \delta^D \circ \varphi \circ \delta$
which proves that $\delta$ is compatible with the natural
commutator pairings on $H$ and $K \times K^D$.
\end{proof}
\noindent
Now the images of $K$ and $K^D$ under $\delta$ equal the
kernels of the lifts of the Verschiebung $V:A^{(q)} \rightarrow A$
and the relative $q$-Frobenius
$F:A^{(q)} \rightarrow A^{(q^2)}$, respectively.
\newline\indent
First assume that $p>2$.
Combining Theorem \ref{key1}, Theorem \ref{key2}
and Proposition \ref{groth} we get sections
\[
u:K \rightarrow G\left( \big( \pol^{(q)} \big)^{\otimes q} \right)  
\quad \mbox{and} \quad
v:K^D \rightarrow G\left( \big( \pol^{(q)} \big)^{\otimes q} \right)
\]
of the natural projection
$G\left( \big( \pol^{(q)} \big)^{\otimes q} \right) \rightarrow H$.
Here $K$ and $K^D$ are considered as subgroups of $H$ via
the Lagrangian structure $\delta$ constructed above.
By Proposition \ref{newtheta} the triple
$(\delta,u,v)$ gives a theta structure of type $K$ for the
pair $\left( A, \big( \pol^{(q)} \big)^{\otimes q} \right)$.
\newline\indent
The above proof applies to the case $p=2$ with some minor
change. Assume that $p=2$.
We claim that there exists a canonical theta structure of type
$(\xz / q \xz)^g_R$ for the pair
\[
\left( A^{(2q)} , \big( \pol^{(2q)} \big)^{\otimes q} \right).
\]
We can argue as above replacing $A$ by $A^{(2)}$.
Theorem \ref{key2}.2 requires the choice of an isomorphism
\[
A^{(2)}[2] \stackrel{\sim}{\rightarrow}
A^{(2)}[2]^{\mathrm{loc}} \times A^{(2)}[2]^{\mathrm{et}}.
\]
We claim that there is a canonical choice.
It is given by the section
\[
A^{(q)}[q]^{\mathrm{et}}
\rightarrow A^{(q)}[q]
\]
induced by the isomorphism
$F[q]^{\mathrm{et}}$ (see diagram (\ref{imp})).
This finishes the proof of Theorem \ref{excantheta}.

\subsection{Proof of Corollary \ref{canlifttheta}}
\label{canliftproof}

We use the notation of Section \ref{thecanonicalthetastructure}.
Assume that $k$ is perfect and $A$ is the canonical lift of $A_k$.
Let $\sigma$ denote an automorphism of $R$ lifting the $(\delta q)$-th power
automorphism of $k$ where $\delta$ is as in Theorem \ref{excantheta}.
We denote by $A^{(\sigma)}$ the pull back of $A$ by the
automorphism $\sigma^{-1}$.
On $A^{(\sigma)}$ there exists an ample symmetric line bundle
$\pol^{(\sigma)}$ of degree $1$ which is defined to be the
pull back of $\pol$ along the projection $A^{(\sigma)}
\rightarrow A$.
\newline\indent
We set $A'=\big( A^{(\sigma)} \big)^{(\delta q)}$.
Since $A'$ is a canonical lift of $A_k$
it follows by uniqueness that there
exists an isomorphism
$\tau:A \stackrel{\sim}{\rightarrow} A'$.
We set
\[
\mathcal{M}=\tau^* \left( \big( \pol^{(\sigma)} \big)^{( \delta q)} \right).
\]
We claim that $\mathcal{M}^{\otimes q } \cong \pol ^{\otimes q}$.
We set $\pol'=\pol \otimes \mathcal{M}^{-1}$.
Our claim follows from the fact that $\mathcal{M}$ and $\pol$ are symmetric
and hence
\[
\langle \pol' \rangle \in \mathrm{Pic}^0_{A/R}[2](R).
\]
By Theorem \ref{excantheta} the line bundle
$\mathcal{M}^{\otimes q }$ has a canonical theta structure
of type $(\xz / q \xz)^g_R$. The latter gives a canonical theta
structure for $\pol^{\otimes q}$.
This completes the proof of the corollary.

\section*{Acknowledgements}

I owe thanks to Bas Edixhoven for assisting me with the proof of
Theorem \ref{excantheta}. Without him this article would
probably not exist.
I'm indebted to Ben Moonen and Michel Raynaud who independently
from each other communicated to me the version of the
proof of Theorem \ref{key2}.1
that we present in Section \ref{second}. 
Their proof supersedes the author's original proof based on an idea
which was the outcome of a discussion with Frans Oort.

\bibliographystyle{plain}

\end{document}